\documentclass[5p,twocolumn]{elsarticle}

\usepackage{lineno,hyperref}
\usepackage{graphicx}
\usepackage{amsmath}
\usepackage{booktabs}
\usepackage{caption}
\usepackage{etoolbox}
\usepackage{subfig}
\modulolinenumbers[5]

\journal{Global Food Security}

\makeatletter
\patchcmd{\ps@pprintTitle}{\footnotesize\itshape
Preprint submitted to \ifx\@journal\@empty Elsevier
\else\@journal\fi\hfill\today}{\relax}{}{}
\makeatother

\begin{document}

\begin{frontmatter}

\title{An independent evaluation of the Famine Early Warning Systems Network food security projections}

\author[myprimaryaddress]{Marco Bertetti\corref{cor1}}
\author[]{Paolo Agnolucci}
\author[]{Alvaro Calzadilla}
\author[]{Licia Capra}
\ead{marco.bertetti.21@ucl.ac.uk}

\address[myprimaryaddress]{Bartlett School of Environment, Energy and Resources, University College London, UK}

\begin{abstract}
Reports from the Famine Early Warning Systems Network (FEWSNET) serve as the benchmark for food security predictions which is crucial for stakeholders in planning interventions and support people in need. 
This paper assesses the predictive accuracy of FEWSNET's food security forecasting, by comparing predictions to the following ground truth assessments at the administrative boundaries-livelihood level, revealing an overall high accuracy of 78\% across diverse timeframes and locations.  
However, our analysis also shows significant variations in performance across distinct regions and prediction periods. Therefore, our analysis sheds light on strengths, weaknesses, and areas for improvement in the context of food security predictions.  The insights derived from this study not only enhance our understanding of FEWSNET's capabilities but also emphasize the importance of continuous refinement in forecasting methodologies.
\end{abstract}


\end{frontmatter}

\section{Introduction}
\label{sec:introduction}
The UN's 2030 Agenda for Sustainable Development is a plan of action that was adopted by all member states of the United Nations in 2015.
One of the agendas's Sustainable Development Goals aims to \textit{ensure access by all people, to nutritious and sufficient food all year round}\cite{UNSDG_Hunger}, requiring  food security to be guaranteed for all. Unfortunately, due to increasing extreme weather patterns, conflicts, economic shocks, and the impact of Covid-19, the number of people suffering from hunger is rising, and the world is not currently on track to end hunger by 2030 (WFP, 2020)\cite{WFP-2020}. By the end of 2022, the FAO's Global Report on Food Crises documented unprecedented levels of acute food insecurity, impacting 250 million individuals, and hunger world-wide has remained at much higher levels compared to pre-COVID-19 pandemic levels (FAO, 2023)\cite{FAO-etal-2023}. \\  

The Famine Early Warning Systems Network (FEWSNET), established by
the United States Agency for International Development (USAID) in 1985 with the specific goal to fight hunger\cite{USAID}, is today considered the gold standard of food security (FS) assessments and predictions. These are generated by applying a comprehensive methodology, which includes data analysis based on heterogeneous data sources, from food prices to weather and rainfall patterns, to socioeconomic indicators, and the application of a series of assumptions, on top of which most likely scenarios can develop\cite{FEWSNET_ScenarioDevelopment}. Reports are compiled and released every four months, and are routinely used by governments and international institutions to plan for aid and prevent the escalation of food crises.
Understanding the accuracy of these projections is essential for policymakers and researchers alike. Policymakers rely on these reports to allocate resources and plan for humanitarian interventions. Researchers use this data to train and evaluate predictive models since at least 2019\cite{Perez2019}\cite{Andree2020}, the reliability of which directly impacts the quality of research and the development of effective solutions.\\

From an academic standpoint, efforts related to evaluating the performance of FEWSNET have been limited. Krishnamurthy et al. (2019)\cite{Krishnamurthy2019} compared projections against later ground truth assessments created by FEWSNET (FEWSNET data explained in the Method section below), using Ethiopia as a case study, showing that FEWSNET's predictions are mostly accurate. However, the study also highlighted variations in performance between regions, as well as a lack of improvement in predictive performance over time. They later expanded the geographical coverage of their assessment to the whole Horn of Africa (Krishnamurthy et al. (2020)\cite{KRISHNAMURTHY2020100374}), showing that the vast majority of the prediction error could be explained by the inability to predict climate-related disasters or the occurrence of conflict and civil unrest. While these studies provided the first independent evaluation of FEWSNET performance, their limited geographical coverage does not allow to draw conclusions on the overall ability by FEWSNET to forecast FS.

Backer and Billing (2021)\cite{BACKER2021100510} further expanded the coverage of FEWSNET's independent evaluation by analysing performance of the projections in 25 countries between July 2009 and June 2020, by comparing the mid-term projections created by FEWSNET against the current assessment from the following report publication period at a 0.5 degree raster grid, therefore creating a stable way of comparing performance over space and time. This work shows that overall accuracy is good, at 84\% but it declines sharply in presence of severe food crises. This paper greatly expanded the geographical coverage of FEWSNET's performance evaluation, and its raster methodology represents a reproducible framework for future comparisons. However, the raster method used in this work diverges from the methodology that FEWSNET uses to generate its predictions and assessments (which is at the administrative boundaries-livelihood level), introducing challenges in understanding the causes behind misclassifications and errors. Moreover, the use of raster data to delineate geographical areas lacks the representation of natural subdivisions, such as social, geographical, or political boundaries. Raster, being an artificial method of territory subdivision, falls short in capturing the complex dynamics of phenomena such as the spread of conflicts or the progression of droughts.\\

This paper examines the accuracy of FEWSNET's predictions across all African countries where the organization operates, specifically focusing on the administrative boundaries-livelihood level, which is the level of aggregation utilized by FEWSNET's methodology when generating their reports, thus providing a directly comparable framework. Building on the prior studies presented above, our work employs a framework that directly aligns with FEWSNET's prediction and food security (FS) assessment methodologies. Following a similar approach to Backer and Billing (2021), we assess FEWSNET's near-term projections by comparing them to the current assessments for subsequent publication periods. This comparative analysis involves calculating various classification metrics (accuracy, precision, recall, F1) for each area and timeframe. Additionally, we enhance our understanding by comparing FEWSNET's projections with simple rule-based models, offering insights into the strengths and weaknesses of the current methodology.\\

This study is the first to assess the accuracy of FEWS NET's assessments and projections using the same level of spatial granularity employed by the network across 24 African countries and Yemen, thereby addressing a significant gap in the literature.

\section{Methods}
\label{sec:methods}
\subsection{Data}
FEWSNET's primary contribution to addressing hunger lies in the regular publication of Food Security Outlook reports. Since 2016, these reports have been issued three times per year, specifically in February, June, and October. Reports are produced for each of the 24 African countries and Yemen where FEWSNET operates, and provide comprehensive insights into the current food security situation, its underlying causes, and future projections. Additionally, the reports outline assumptions about how these conditions will evolve and their anticipated impact on food security in the coming six months. \\

\begin{figure*}[h]
\centering
\includegraphics[width=0.5\textwidth]{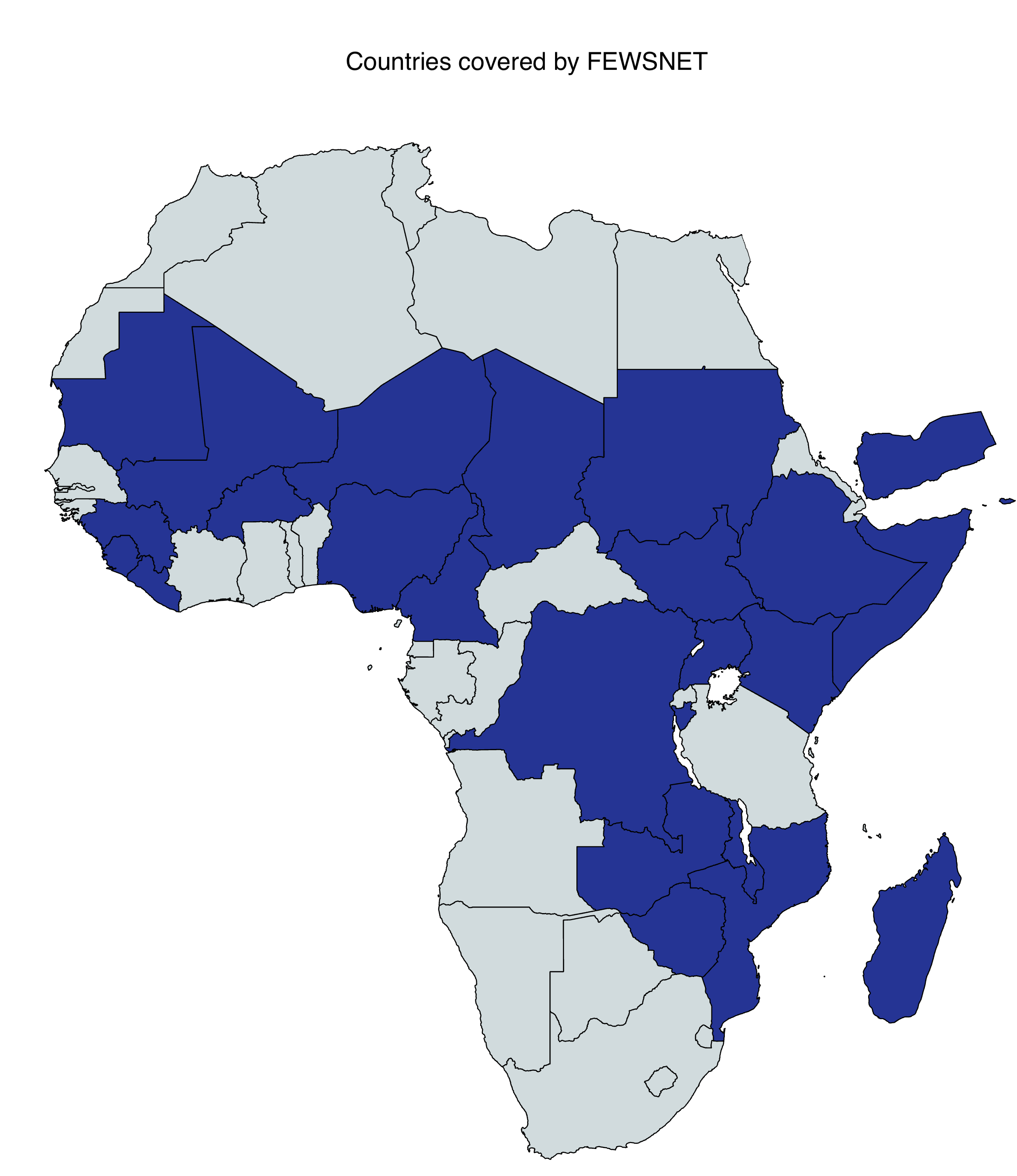}
\captionsetup{labelfont=bf}
\caption{Countries in blue are those for which FEWSNET generates regular assessments and projections.}
\end{figure*}

FEWSNET employs a comprehensive methodology in the creation of its reports, a key aspect of which involves the assessment of FS in various geographical areas being surveyed. This assessment utilises the Integrated Food Security Phase Classification (IPC) scale, which is a five-class system designed to provide a detailed understanding of food security levels. The IPC scale spans from Class 1, signifying ``minimal'' food insecurity, to Class 5, denoting the most severe condition, ``Famine''. Throughout this analysis, the terms ``phase'' and ``class'' will be used interchangeably to refer to the different levels on the IPC scale. The table below provides a detailed description of each phase.

\begin{figure*}[h]
\centering
\includegraphics[width=0.8\textwidth]{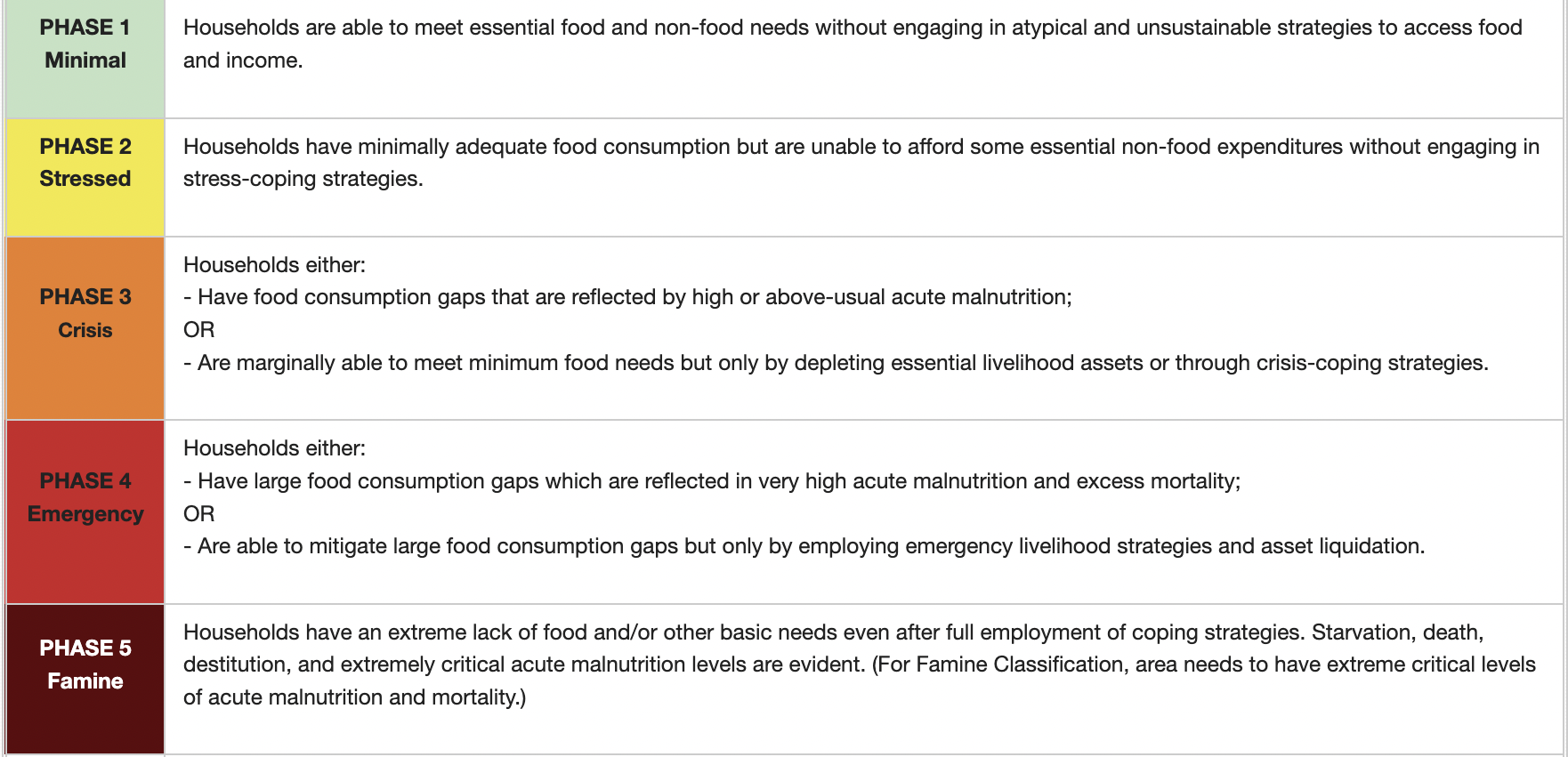}
\caption{Source: fews.net}
\end{figure*}

With each publication cycle, FEWSNET also releases the raw data behind each report, enabling researchers to conduct independent work on the FS topic. The key datasets made available by FEWSNET, and those we have utilised to create our baselines are the following:

\begin{itemize}
  \item \textbf{Food Security Classification}: Provides the IPC phase classification for each geographical area. Specifically, this data includes:
  
  \begin{itemize}
  \item \textbf{CS}: The Current Situation assessment for the state of FS at the time of the reports' publication;
  \item \textbf{ML1}: Short-term projection for the state of FS for the three months following the reports' publication;
  \item \textbf{ML2}: Medium-term projection for the state of FS for the six months following the reports' publication; 
  \end{itemize} 
    
Each file contains one row for each of the five IPC FS classes. Associated with each class is a set of polygons that delineate specific regions on the map. These polygons define the boundaries of areas classified under each FS state, enabling the visualization of various areas according to their food security status.
  
  \item \textbf{Administrative Boundaries}: shapefiles containing polygons with information about sub-national level administrative boundaries; 
  \item \textbf{Livelihood Zones}: shapefiles containing polygons that define geographical areas of each country where people share similar means of obtaining food sources, either by direct production or by income-generating activities, such as ``High Food Production and Surplus Trade'' or ``Agropastoral Zone''; \\

  The Administrative Boundaries and Livelihood Zones shapefiles serve different purposes and represent different types of geographic information. While Administrative Boundaries focus on political and administrative divisions within a country, Livelihood Zones focus on geographic areas with similar livelihood characteristics. In practice, the areas in these files can intersect and overlap; for example, a specific region may have a coastal area where the main livelihood is fishing, and and inland area where crops are cultivated, or multiple areas in different countries can share the same livelihood type (i.e. trade).

\end{itemize} 

While FEWSNET has been providing data since its inception in 2009, this study specifically focuses on reports published from 2016 onwards, for three key reasons. First, in 2016 FEWSNET transitioned from a quarterly publication schedule to the current practice of releasing reports three times per year. Only utilising data since 2016 ensures consistency in the comparison of reports over time. Secondly, as the purpose of this work is to assess the accuracy of FEWSNET predictions, which are being used as a benchmark against which researchers compare their models, focusing on the most recent years provides the most relevant time frame against which to build a baseline. \\
Lastly, as FEWSNET continually refines its methodologies to improve the accuracy and effectiveness of its FS predictions and assessments, and the quality and availability of data pertinent to these predictions have significantly increased in recent years, focusing on the latest data establishes the most challenging performance baseline. 

\subsection{Methodology}
To evaluate the accuracy of FEWSNET projections, a robust framework is essential, allowing for the calculation of diverse metrics. To establish this framework, we merged the data outlined in the preceding section, forming a geospatially consistent dataset. This dataset serves as the foundation for analyzing the performance of FEWSNET predictions across both temporal and geographical dimensions. In order to ensure a consistent spatial aggregation over time, we compared the near-term projection (ML1) published in each report period with the subsequent report's current assessment (CS). This approach provides a dynamic perspective, allowing for a nuanced evaluation of FEWSNET's predictive capabilities. To achieve spatial consistency, we mapped the spatial data from each report onto the latest administrative boundaries and livelihood zones pairs made available by FEWSNET. These boundaries, known for their inherent stability, only change over extended periods, ensuring a reliable spatial area for comparison. The heuristics-based models we created, detailed in a later section, serve as a comparison against FEWSNET's predictions and follow the same methodology for the metrics' calculation.

\begin{figure*}[h]
\centering
\includegraphics[width=0.9\textwidth]{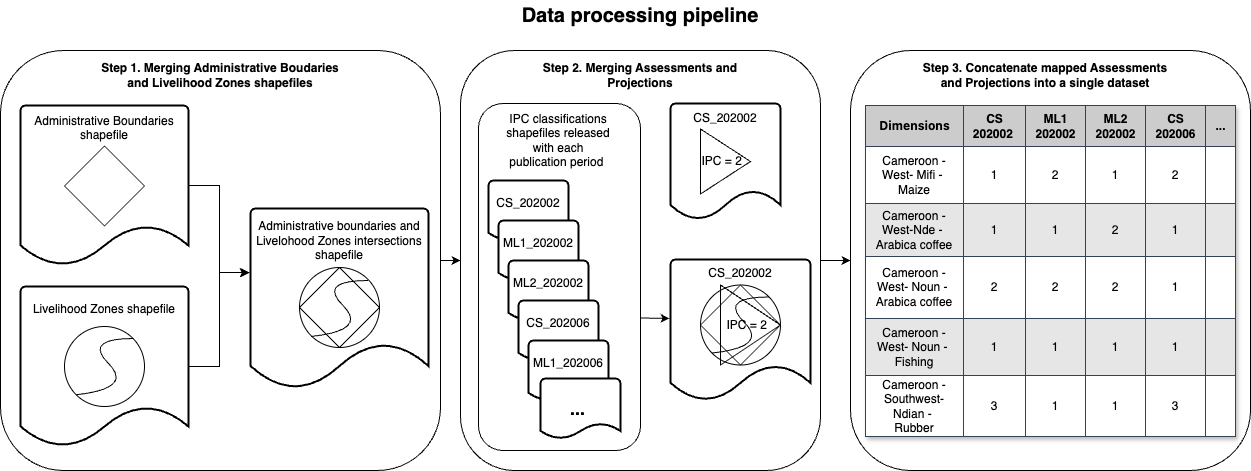}
\caption{Diagram of the data processing pipeline: in Step 1 we combine Administrative Boundaries and Livelihood Zones, in Step 2 we merge the output from the previous step with each published assessment and projections; Step 3 combines all the data into a single table, covering assessments and projections over time and geographies}
\end{figure*}

The following subsections provides a detailed explanation of the computational steps we have taken to process the data, as well as the logic used to create the rule-based models.

\subsubsection{Data processing}
The first step combines the data from the administrative boundaries and the livelihood zones. This was accomplished by intersecting the polygons from each set of data, creating more, smaller areas containing information from both sources (Fig. 3, Step 1).

Since both administrative boundaries and livelihood zones seldom change, their combination provides a stable and granular geo-spacial dimension for comparisons over time; we have used the most recent files (as of September 2023) as fixed point for the mapping. \\

The figures below show and example of the raw data contained in each dataset. 

\begin{figure*}[h]
  \centering
  \subfloat{\includegraphics[width=0.4\textwidth]{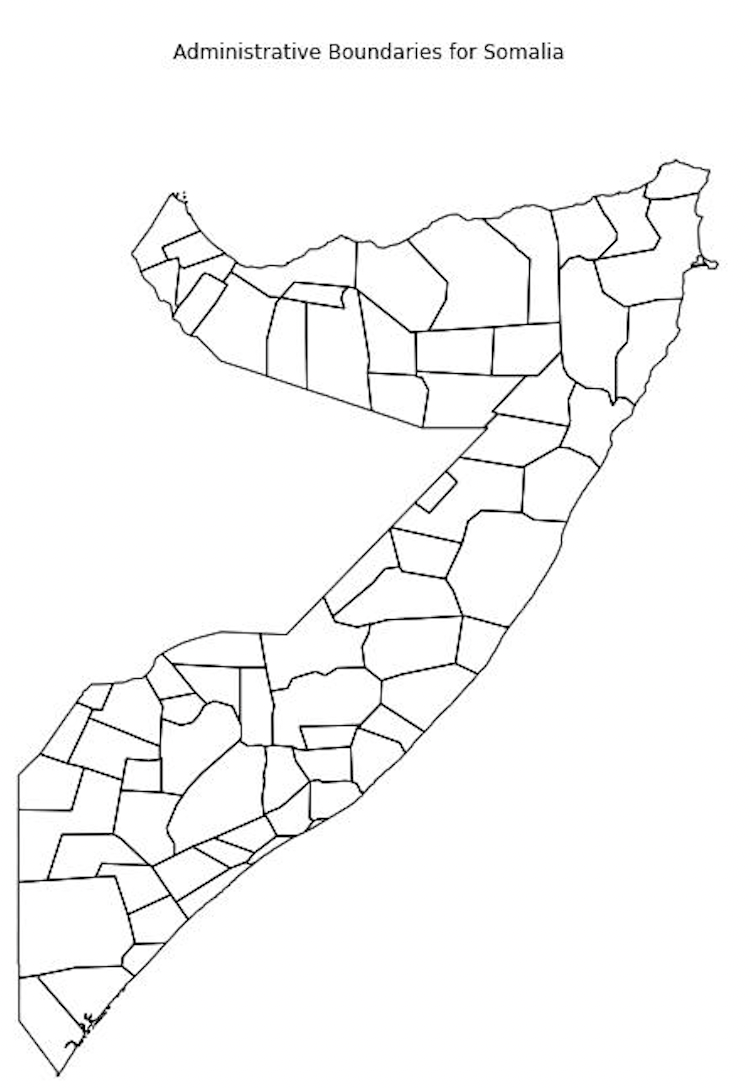}}
  \hfill
  \subfloat{\includegraphics[width=0.4\textwidth]{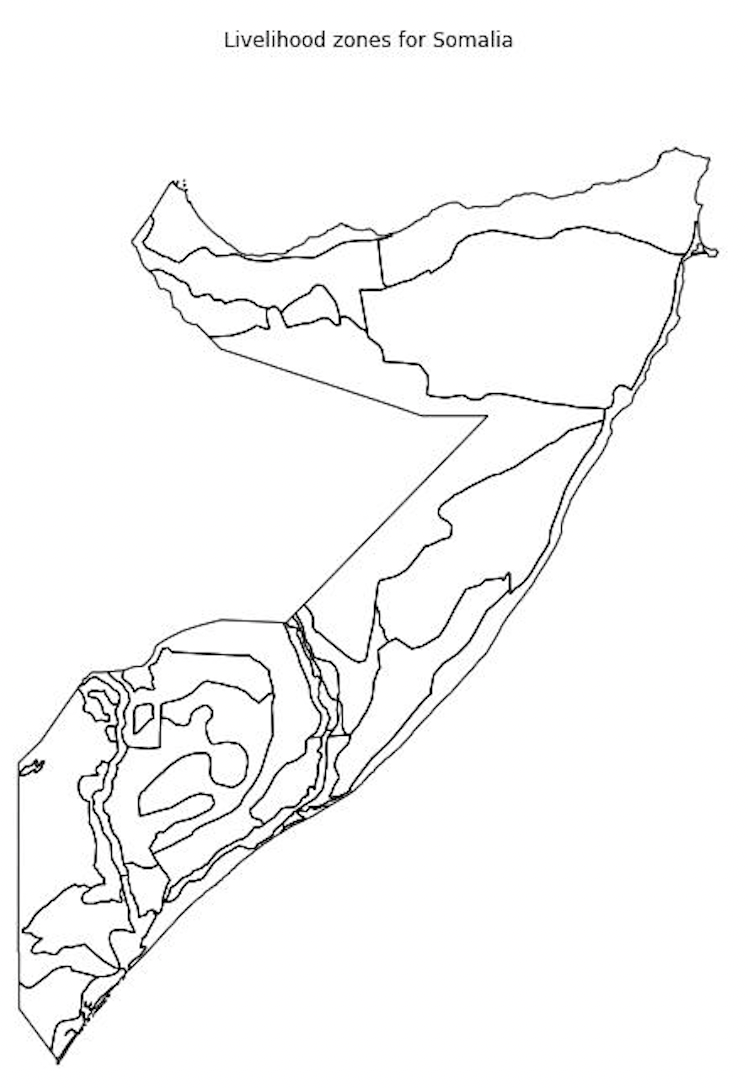}}
  \caption{Example of visualisation of Administrative Boundaries and Livelihood zones in Somalia from FEWSNET shape files.}
\end{figure*}

The second steps consists in iterating across each of the Food Security Classification files. Since this data also consists of shapefiles, we could again intersect each of the polygons containing the IPC classification for each period of publication with the polygons containing administrative boundaries and livelihood zones information (Fig. 3, Step 2)\footnote{Total number of merges is 190. This is equal to 7 years of data * 3 publications per year * 3 files per publication (CS, ML1, ML2) * 3 geographical areas (SA, EA, WA) + the initial merge between livelihood zones and administrative boundaries.}.

After computing the intersection, we have removed polygons with an area smaller than 0.005; these very small polygons can occur due to the non-perfect overlapping between the original shapes, and can cause duplication in the FS-Admin-Livelihood combinations; removing such small areas resolved this problem without eliminating any significant amount of land coverage\footnote{This removes less than 0.01\% of the total area.}.\\

The output of this process are new shapefiles, covering the same geographical area as the original Food Security Classification files from FEWSNET, but with an increased number of polygons, each representing the overlap between FS classification, livelihood zones, and administrative boundaries. The figure below shows an example of the output data. \\

\begin{figure}[h]
\centering
\includegraphics[width=0.4\textwidth]{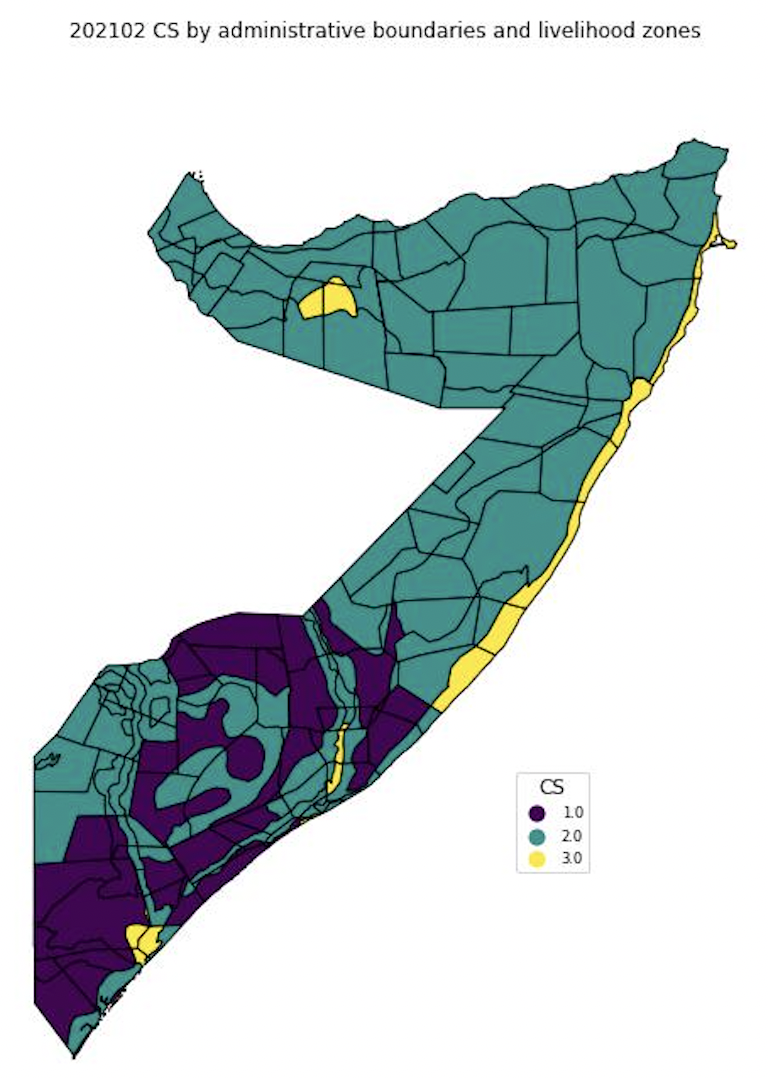}
\caption{Visualisation of the final output dataset, result of combining Administrative Boundaries, Livelihood Zones, and IPC FS assessment shape files published by FEWSNET (CS for February 2021)}
\end{figure}

Lastly, e removed duplicates from the data. We aggregated the data by each Admin-Livelihood combination and, in case any duplicate occurred, we took the highest (worse) FS score for each of these areas (Fig. 3, Step 3)\footnote{This operation affects between 3 and 5\% of the dataset (depending on the period), where the same Admin-Livelihood have different FS scores.}.\\

This method aggregates data based on unique combinations of IPC class, administrative boundaries, and livelihood zones, offering several advantages over raster-based approaches commonly used in the literature (e.g., Backer and Billing, 2021). Firstly, it aligns with FEWSNET's prediction granularity. FEWSNET's FS assessment and prediction involves using a combination of socioeconomic data at a district level, and data on how different crops (livelihoods) have been impacted by most recent weather events, trade and food prices. This approach enables a direct assessment of their food security forecasts, which are made at the district level. 

Secondly, this approach enhances the accuracy of performance metric calculations. FEWSNET's coverage often includes large, single-crop areas, which, if represented using rasters, result in multiple data points. Admin-livelihood aggregation condenses these areas into a single data point, leading to a more fair evaluation.

Lastly, this data aggregation method simplifies data interpretation and communication. By grouping information according to administrative boundaries and livelihood zones, it provides a more intuitive way to convey results and findings to policymakers and stakeholders. 

\subsubsection{FEWSNET predictions and heuristic models}
As previously mentioned, we have evaluated FEWSNET's actual projections and established three distinct rule-based models to serve as baselines for comparison. \\
Below is a detailed description of these baselines and how they relate to the actual projections:

\begin{itemize}

\item \textbf{FEWSNET Projections (FEWSNET)}: The actual projections published by FEWSNET. Assessing the performance of simple models against FEWSNET's actual projections provides a benchmark for evaluating the quality and reliability of alternative forecasting methods. 

\item \textbf{Previous Period's IPC Score (PPS)}: This first model assumes a scenario of no change in food security conditions and simply replicates the Integrated Food Security Phase Classification (IPC) score from the preceding period (ML1t = CSt-1). It provides insight into the stability of food security.

\item \textbf{Same Period Last Year (SPLY)}: Recognizing the influence of seasonal factors on food security, this model considers the IPC score from the same period in the previous year (ML1t = CSt-3). By accounting for seasonality, it offers a valuable perspective on the cyclical nature of food security predictions.

\item \textbf{Maximum of Previous Two Periods (Max-2PP)}: This model selects the highest (worst) IPC score from the previous two periods (ML1t = Max(CSt-1, CSt-2). This approach takes into account recent trends and variations in food security, highlighting potential deteriorations.

\end{itemize} 

The following section will provide a detailed description of the different metrics that will be used to assess the current performance of FEWSNET, and to compare it different prediction methods.

\section{Metrics}
In order to fairly compare the performance of different methods, four main metrics will be used: accuracy, precision, recall, and F1. In this study, all metrics have been computed at the administrative boundaries-livelihood areas level, then aggregated as needed to show national and temporal trends. These are all metrics commonly used to measure the performance of machine learning (ML) classification problems, and they are defined as follows:

\begin{itemize}

\item \textbf{Accuracy} is a measure of how well a classification model correctly predicts the target variable. It is calculated by dividing the number of correct predictions by the total number of predictions. Accuracy provides a general overview of the model's correctness. However, it can be misleading in imbalanced datasets, as high accuracy can still be achieved by predicting the majority class while neglecting the minority class. This is very much the case when working with FS where, thankfully, crisis situations are the vast minority of total cases.

\item \textbf{Precision} focuses on the positive predictions (when the model predicts a positive class label) made by the model and measures the ability to avoid false positive predictions. It is calculated as the ratio of true positives to the observations predicted as positives. Precision is particularly useful in situations where false positives are costly or undesired. 

\item \textbf{Recall} also known as sensitivity or true positive rate (when the model correctly predicts the positive class), measures the ability of a model to correctly identify positive instances from the total number of actual positive instances. It is calculated as the ratio of true positives to the observations that are actually positive. Recall is valuable in scenarios where false negatives are critical. For example, in FS, it is key to not miss crisis situations.

\item \textbf{F1 score} is the harmonic mean of precision and recall. The harmonic mean gives equal weight to precision and recall, making the F1 score useful in situations where both metrics are important. The F1 score ranges from 0 to 1, with 1 representing perfect precision and recall, and 0 indicating poor performance in either metric.
\end{itemize} 

In the next section, we describe the results we have obtained from analysing FEWSNET data, providing a full assessment of their predictive performance.

\section{Results}
\label{sec:results}
Our analysis shows that FEWSNET's predictive accuracy is high, with an average accuracy of around 78\% across time series and geographies, aligning with findings from related studies (Backer and Billing, 2021). For instance, Figure 6 illustrates the distribution of these predictions. Most errors were minor, with 98.7\% within ±1 IPC class. Only 1.3\% of errors deviated by two or more classes. Specific challenges in areas like Middle Shabelle, Somalia, and Tahoua, Niger, illustrate the impact of local adversities. In Middle Shabelle, agricultural setbacks such as drought, pest damage, and financial constraints significantly lowered crop yields, especially in the Cowpea Belt Agropastoral livelihood zone. In Tahoua, intensified terrorist activities disrupted humanitarian efforts, complicating access to displaced populations and restricting intervention scopes due to security measures.

\begin{figure}[h]
\centering
\includegraphics[width=0.5\textwidth]{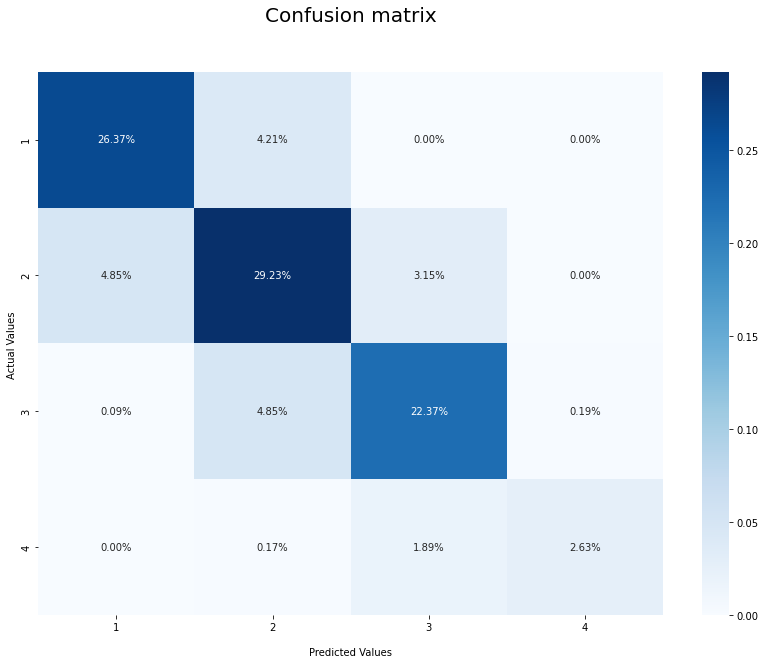}
\caption{Confusion matrix showing the accuracy for the predictions made by FEWSNET in February 2022 vs the actuals computed in June 2022. Most predictions are accurate, and the vast majority of the rror falls within +-1 predicted class.}
\end{figure}

Focusing on the temporal dimension, we  observe that FEWSNET's predictive performance exhibits periodic fluctuations, with a slight decline averaging 3\% over recent years (Fig 7) specifically when comparing the period between the 2017-2019 and 2020 onwards. This is particularly interesting, especially considering the increase in quality, availability, and frequency of data over the past few years.\\

When we compare FEWSNET's predictions against our rule-based models over time, as we would expect FEWSNET's outperform all the simple heuristics-based models, whose accuracy averages at 0.76\% for PPS, for 74\% for Max-2PP and 67\% for SPLY. Additionally, performance of rule-based models show greater variabililty across periods than FEWSNET's with minimum levels of accuracy as low as 53\% for SPLY, and 64\% for Max-2PP. 

\begin{figure*}[htb]
  \centering
  \includegraphics[scale=0.4]{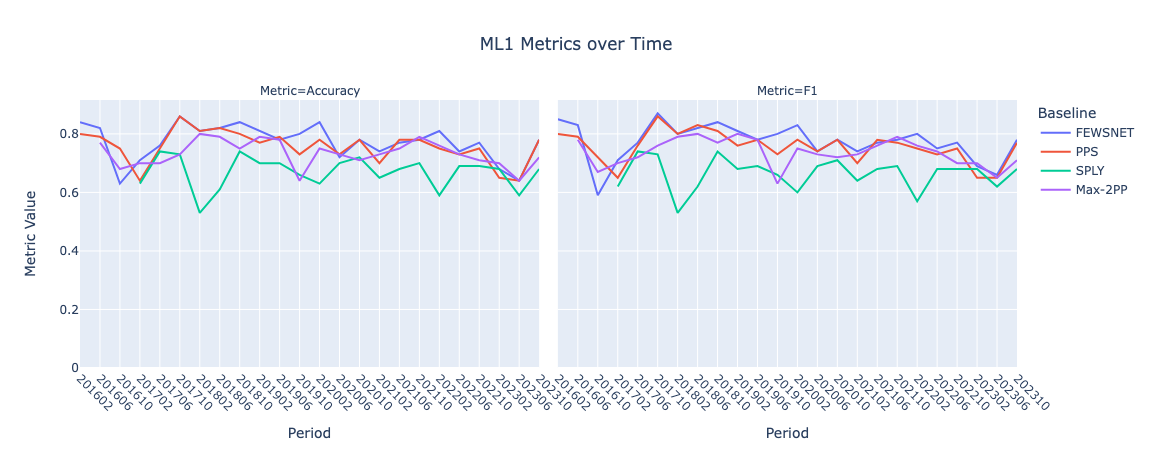}

  \vspace{2ex} 

  \includegraphics[scale=0.4]{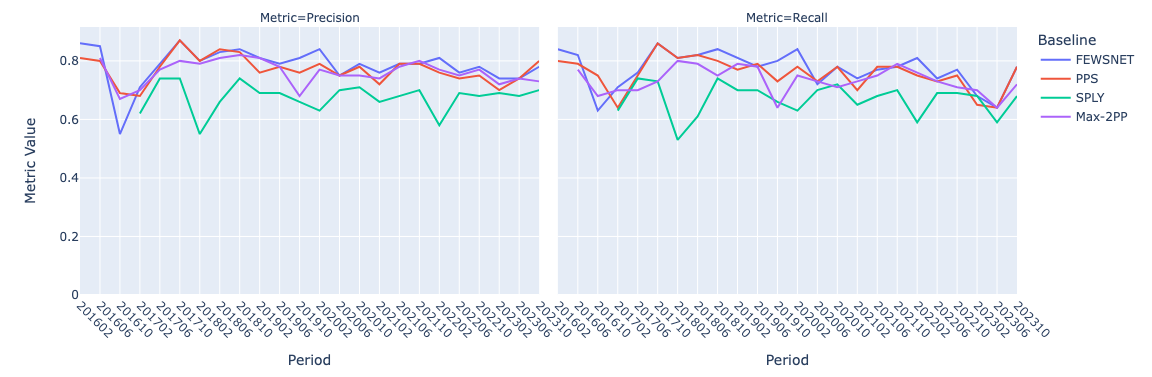}
  \caption{Accuracy, Precision, Recall, and F1 over time across all administrative boundaries-livelihood areas covered by FEWSNET assessments and projections. We can see periodic fluctuations and an overall slight decrease in performance since 2020.}
\end{figure*}

Since aggregating data by period can hide local variations, we analysed the performance of predictions over different geographies, noticing great variations between different countries (Fig 8); Somalia and South Sudan show the worse performance, with overall average accuracy of 57\% and 60\% respectively. For example, in the second half of 2019 the prediction accuracy in Somalia was amongst the lowest ever recorded, due to widely inaccurate rainfall pattern predictions where, after a prolonged period of draughts, rainfall occurred at a much above average rate causing flash floods and a much higher than expected rated of displacement. In the case of South Sudan, most of the consistently low prediction accuracy can be attributed to the occurrence of conflicts. Violence not only erupts quickly and unexpectedly and is therefore challenging to include into FS predictions, but is also poorly explored and understood within the FS space, as highlighted by several scholars within the literature (Brown et al. (2020)\cite{Brown2020}, (Brown et al. (2021)\cite{brown2021}, Kemmerling et al. (2022)\cite{kemmerling2022}).\\

In countries such as Uganda and Cameroon, on the other hand, predictions are more accurate, at around 92\%. This is due to the overall higher stability of FS in these countries. For example, Uganda is food secure most of the time therefore predictions of IPC = 1 will be correct most of the time. In Cameroon, where FS is unfortunately not as high, the drivers behind FS are more stable and better understood, such as the ongoing below average production and the loss of subsistence methods experienced by people in the Northwest and Southwest regions of the country, therefore making predictions more accurate. \\

It is also worth noting that there are gaps in the predictions and assessment generated by FEWSNET; some countries such as Sierra Leone and Zambia do not have any data available since 2016 and 2017 respectively, while for other countries assessments were first generated only in 2020.

\begin{figure*}[h]
\centering
\includegraphics[width=\textwidth]{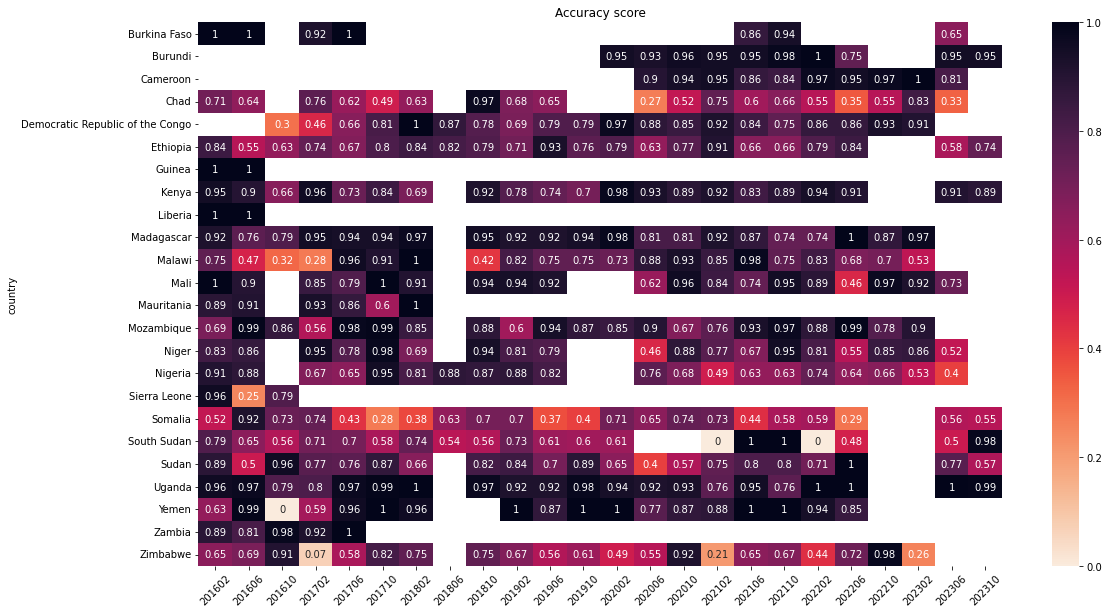}
\caption{Accuracy of FEWSNET projections by country and period. We can observe large variations in performance both between countries and across time.}
\end{figure*}

FEWSNET's predictive performance, while overall fairly accurate, show variations both over time and across countries.
Several factors may contribute to this. Firstly, external events, such as the global COVID-19 pandemic, have introduced unprecedented challenges, disrupting established patterns in food security. The pandemic's economic repercussions, supply chain disruptions, and changes in consumer behavior could have impacted the reliability of the traditional FEWSNET's methodology.\\

Looking at performance over time, the greatest drop in performance is observed in the predictions generated in the summer of 2020, predicting the state of FS for the following fall. In this particular case, most of the error was due to an overestimation of the severity of the near-term impact of COVID-19 on FS, especially in countries such as Yemen and Sudan, where an ``Emergency'' level of FS was predicted but, thankfully, did not occur. \\

At a spatial level, overall performance may improve by collecting more country-specific, near-real-time data, in order to better include, for instance, the insurgence of sudden conflicts and their impact, and to better understand the country's adaptability and resilience to extreme weather events.

\section{Discussion}
\label{sec:discussion}
This paper introduces a robust and replicable approach to evaluating the performance of FEWSNET's projections. The novelty and key contribution of our work lie in the spatial granularity of the framework, specifically at the administrative boundaries and livelihood zones level. This level aligns with FEWSNET's own methodology for creating assessments and projections, establishing a direct link to FEWSNET's reports. Our analysis reveals an overall accuracy of 78\% in FS projections across all areas. Additionally, the study brings attention to temporal and spatial variations in performance, providing valuable insights for future improvements. \\

This research offers substantial potential impact in both governmental and academic realms. In government and policy contexts, more focus and resources can be allocated to areas where variability in FS predictive performance is higher, for instance by a targeted increase in the assessment frequencies, or by integrating more data related to key vulnerabilities, such as conflicts. Such informed decision-making can enhance the efficiency and effectiveness of interventions, ensuring that efforts are concentrated where they are most needed, thereby contributing to more robust food security policies.
On the academic front, this work improves the understanding of FEWSNET's performance, today's FS predictions gold standard, used by many scholars as a benchmark for their own predictive models. The vulnerabilities highlighted can also be interpreted as an ``error analysis'' on the current FS projections, and can help researchers selecting data for their models, improving prediction accuracy over time. We also call for more research on the analysis of the impact of conflicts on the state of FS, as well as for more data collection covering topics such as infrastructure resilience, a key aspect in better predicting the state of FS in case of extreme events.
\\

While this work provides valuable insights, it is essential to underscore some limitations inherent in the current approach. Firstly, the methodology relies entirely on data published by FEWSNET, making it susceptible to any changes in data format, methodology, or assumptions adopted by FEWSNET. Therefore, adjustments in data processing may be necessary to accommodate any alterations, should they occur. Secondly, the approach assigns identical weight to all geographical areas, irrespective of variations in size and population density. There exists an opportunity for additional refinement, such as shifting from accuracy defined as percentage of correctly classified areas to percentage of correctly classified households. Lastly, this analysis could be extended  beyond Africa to encompass all countries assessed by FEWSNET.
\\

It is also important to note that alternative approaches could be adopted to achieve the same goal. For example, one could compare FEWSNET's predictions against correlated metrics collected by different entities, such as the Food Consumption Score (FCS) or the Coping Strategies Index (CSI) collected by the WFP. These metrics are also being used as target variables for predictive models\cite{Lentz2019}\cite{Zhou2021}\cite{Deléglise2020}\cite{Deléglise2021}\cite{WFPVAM}. However, individual metrics do not provide an accurate representation of the overall state of FS, but merely cover one specific aspects of its multifaceted nature. Moreover, many of these metrics lack a predictive element, making it challenging to assess their forecasting accuracy, a crucial factor in addressing long-term FS.

Evaluating FEWSNET's predictions against economic indicators related to food security, such as local market prices, income levels, or agricultural production data, can also validate its forecasts from an economic standpoint. However, comparing FS to economic indicators may oversimplify the multifaceted nature of food security, reducing it to a singular economic perspective, which is likely to fall short of capturing the intricate interplay of social, environmental, and political factors that contribute to the overall resilience and vulnerability of food systems. 

Alternatively, conducting case studies or field validations in specific regions could verify FEWSNET's predictions against ground-truth data, providing a detailed understanding of its performance in varying contexts. However, this entails substantial costs, both in terms of finances and time. Conducting comprehensive studies in multiple locations can be costly and time consuming. This restricts their scalability and feasibility for widespread application, limiting their utility in evaluating FS across diverse regions or conducting repeated assessments over time.

\section{Conclusion}
\label{sec:conclusion}
This paper underscores the critical importance of independent evaluations of FEWSNET's food security predictions, offering a valuable contribution to the ongoing discourse surrounding global hunger and sustainable development. The study's in-depth analysis at the administrative boundaries - livelihood level reveals a nuanced perspective on the strengths and challenges of FEWSNET's methodology.

The findings emphasize the need for continuous scrutiny and validation of food security prediction models, especially as the world faces increasingly complex challenges such as the COVID-19 pandemic, conflicts, and extreme weather events. While the results presented in this paper indicate that FEWSNET's predictions generally perform well, with a fairly high level of accuracy, the slight decline in predictive performance observed in recent years highlights the dynamic nature of the factors influencing food security and underscores the necessity for more adaptive and resilient prediction methodologies.

Moreover, the study advocates for additional investigation into the impact of conflicts on food security predictions. The unpredictability and rapid escalation of conflicts pose significant challenges to accurate forecasting, as demonstrated by the varying performance across different regions, with conflicts being a major contributing factor to lower accuracy levels. Understanding the intricate dynamics between conflicts and food security is crucial for refining prediction models and developing targeted interventions in conflict-prone regions.

\section*{Declaration of generative AI and AI-assisted technologies in the writing process}

During the preparation of this work the author used ChatGPT in order to correct typos and improve readability of the article. After using this tool, the author reviewed and edited the content as needed and takes full responsibility for the content of the publication.

\bibliography{sources.bib}

\end{document}